\documentstyle[12pt]{article}
\newenvironment{beweis}{{\it Proof.}\ }{\ $\ \ \ \Diamond$ \\\ }

 \newcounter{nsatz}[section]
 \newcounter{nlemma}[section]
 \newcounter{ndef}[section]
 \newcounter{nhyp}[section]
 \newcounter{nconjecture}[section]
 \newcounter{ncor}[section]
 \newcounter{nrem}[section]
 \newcounter{nexample}[section]
 \newcounter{nprop}[section]

 \newenvironment{nsatz}{\refstepcounter{nsatz}{\bf \arabic{section}.\arabic{nsatz}}\
                {\sc\bf Theorem.\ }\it}{\\\\ \rm}

 \newenvironment{nlemma}{\setcounter{nlemma}{\value{nsatz}}
                \refstepcounter{nlemma}
                \setcounter{nsatz}{\value{nlemma}}
                {\bf \arabic{section}.\arabic{nsatz}}\
                {\sc\bf Lemma.\ }\it}{\\\\ \rm}

 \newenvironment{ndef}{\setcounter{ndef}{\value{nsatz}}\refstepcounter{ndef}
                \setcounter{nsatz}{\value{ndef}}
                {\bf \arabic{section}.\arabic{nsatz}}\
                {\sc\bf Definition.\ }}{\\\\ \rm}

 \newenvironment{ncor}{\setcounter{ncor}{\value{nsatz}}
                \refstepcounter{ncor}
                \setcounter{nsatz}{\value{ncor}}
                {\bf \arabic{section}.\arabic{nsatz}}\
                {\sc\bf Corollary.\ }\it}{\\\\ \rm}

\begin{document}
\newcommand{\n}{{\mbox{\rm I$\!$N}}}
\newcommand{\z}{{\mbox{{\sf Z\hspace{-0.4em}Z}}}}
\newcommand{\R}{{\mbox{\rm I$\!$R}}}
\newcommand{\Q}{{\mbox{\rm I$\!\!\!$Q}}}
\newcommand{\C}{{\mbox{\rm I$\!\!\!$C}}}
\newcommand{\p}{{\mbox{\rm I$\!\!\!$P}}}
\newcommand{\ug}{\ \raisebox{-.3em}{$\stackrel{\scriptstyle \leq}
{\scriptstyle \sim}$} \ }

\thispagestyle{empty}
\setlength{\parindent}{0pt}
\setlength{\parskip}{5pt plus 2pt minus 1pt}

\thispagestyle{empty}
\renewcommand{\thefootnote}{\fnsymbol{footnote}}
\newcommand{\Irr}{{\mbox{\rm Irr}}}
\newcommand{\sIrr}{{\mbox{\scriptsize\rm Irr}}}
\newcommand{\Char}{{\mbox{\rm Char}}}
\newcommand{\Gal}{{\mbox{\rm Gal}}}
\newcommand{\GF}{{\mbox{\rm GF}}}
\newcommand{\Sym}{{\mbox{\rm Sym}}}
\newcommand{\X}{{\mbox{$\setminus$}\mbox{$\!\!\!/$}}}
\newcommand{\dl}{\mbox{\rm dl}}
\mbox{\vspace{4cm}}
\vspace{4cm}
\begin{center}
{\bf \Large\bf  Finite groups have even more conjugacy classes\footnote{This research was supported by the National Security
Agency, Standard Grant No. 08G-206}
}\\
\vspace{3cm}
by\\
\vspace{11pt}
Thomas Michael Keller\\
Department of Mathematics\\
Texas State University\\
601 University Drive\\
San Marcos, TX 78666\\
USA\\
e-mail: keller@txstate.edu\\
\vspace{1cm}
\end{center}
\thispagestyle{empty}
\newpage

\begin{center}
\parbox{12.5cm}{{\small
{\sc Abstract.}
In his paper "Finite groups have many conjugacy classes" (J. London Math. Soc (2) 46 (1992), 239-249),
L. Pyber proved the to date best general lower bounds for the number of conjugacy classes of a finite group
in terms of the order of the group. In this paper we strengthen the main results in Pyber's paper.

}}
\end{center}
\normalsize

\section{Introduction}\label{section0}

The subject of this paper are lower bounds for the number $k(G)$ of
conjugacy classes of a finite group $G$ in terms of the group order.
Over the years this fundamental and venerable problem has attracted
the interest of a number of people, and progress has been slow, but steady, since
Landau's well-known first observation in 1903 that $k(G)$ has to grow with $|G|$.
In 1968, Erd\"{o}s and Tur\'{a}n proved that $k(G)>\log_2\log_2|G|$, and since then
this bound has been improved continuously for various classes of groups. We refer
to E. A. Bertram's excellent survey \cite{bertramsurvey} for more details on the history
of this problem and the results obtained along the way. Suffice it to say here, that it is
conjectured that a bound of the form
\[(*)\quad k(G)>C\log|G|\]
(for a constant $C$) is valid in general, and that there exist $p$-groups with
$k(G)<(\log|G|)^3$.\\
The currently best possible general bounds can be found in L. Pyber's paper
{\it Finite groups have many conjugacy classes} \cite{pyber}, and the goal of this
paper is to improve on the main results in \cite{pyber}.
We prove\\

{\bf Theorem A.} {\it There exists a (explicitly computable) constant $\epsilon_1>0$
such that every finite group $G$ with $|G|\geq 4$ satisfies
\[k(G)\ \geq\ \epsilon_1\frac{\log_2|G|}{(\log_2\log_2|G|)^7},\]
and if in addition $G$ is solvable, then even
\[k(G)\ \geq\ \epsilon_1\frac{\log_2|G|}{\log_2\log_2|G|}.\]}
(See \ref{ncor2.2} and \ref{ncor2.3} below.)\\

Pyber obtained the same results with an exponent of 8 instead of 7 in the first part,
and with a denominator of $(\log_2\log_2|G|)^3$ instead of $\log_2\log_2|G|$ in the second
part.\\

Two remarks are in order here. First, since our proof of Theorem A largely follows
Pyber, the result for arbitrary groups depends on the Classification of Finite Simple
Groups.\\
Second, as Pyber points out in his paper, a lower bound of the form $\epsilon
\frac{\log|G|}{\log\log|G|}$ is the best possible that can be achieved with his
approach. Thus our lower bound for solvable groups fully exhausts Pyber's approach,
and any further progress will most likely have to be based on new stronger lower bounds for nilpotent
groups than currently available.\\

Theorem A will be a consequence of the following improvement of Pyber's Theorem B in
\cite{pyber}:\\

{\bf Theorem B.} {\it There is a (explicitly computable) constant $\beta>0$ such that for
any solvable group with trivial Frattini subgroup we have
\[k(G)\ >\ |G|^\beta.\] }
(see \ref{nsatz2.1} below.)\\

Theorems A and B will be proved in Section \ref{section2} of the paper. In Section \ref{section3} we will draw some
consequences from Theorem B for solvable groups, thereby generalizing and improving some
results in \cite{bertramsurvey}. Perhaps the most interesting result of these is the
following.\\

{\bf Theorem C.} {\it Let $G$ be a finite solvable group, and write $c$ for the nilpotency
class of the Frattini subgroup of $G$. Then
\[k(G)\ >\
\left(\frac{1}{2}\ c\ |G|^\frac{1}{c}\right)^\frac{\beta}{3},\]
where $\beta$ is the constant from Theorem B.}\\

Notation: All groups in this paper are finite. If $G$ acts on a set
$\Omega$, we write $n(G,\Omega)$ for the number of orbits of $G$ on $V$. If
$V=V(n,q^f)$ is an $n$-dimensional vector space over $\GF(q^f)$ ($q$ a prime power),
then as in \cite{manz/wolf} we write
$\Gamma(q^f)=\Gamma(V)=\{x\mapsto ax^\sigma|a\in\GF(q^m)^x,
\sigma\in \Gal (\GF(q^m)/\GF(q))$. Also $\log$ stand for $\log_2$. For
any group $G$ we write $k(G)$ for the number of its conjugacy classes. If $G$ is a permutation
group on a set $\Omega$, write $s(G)$ for the number of orbits of $G$ on the power set
of $\Omega$. For any group, $F(G)$ and $\Phi(G)$ are its Fitting and Frattini
subgroup, respectively.

\section{Preliminaries}\label{section1}

In this section, for the convenience of the reader, we collect several known results
from the literature that will be needed in the proof of the main
results.\\

We will use a result of S. Seager \cite{seager}, which we slightly reformulate here as
in \cite[Theorem 2.3]{kellersurvey} (note, however, that there is a typo in (ii) in
the statement there).\\

\begin{nsatz}\label{nsatz1.1}
Let $G$ be a solvable group and $V$ a finite faithful irreducible $G$-module.
Then one of the following holds:\\

(i) $|V|\ \leq\ \left(\frac{n(G,V)+1.43}{24^\frac{1}{3}}\right)^c$ for $c=36.435663$,
so in particular, $n(G,V)\geq|V|^\frac{1}{37}$, or\\
(ii) for some integers $m,k$ with $k\leq 0.157\log_3\left(\frac{n(G,V)+1.43}{24^\frac{1}{3}}\right)
\leq\frac{1}{5}\log_3n(G,V)$
and a prime $p$ we have $|V|=p^{mk}$ and $G\ug\Gamma(p^m)\wr S_k$
in its natural action on $V(m,p)^k$.
\end{nsatz}
\begin{beweis}
See Theorem 1 in \cite{seager}.
\end{beweis}\\

\begin{ndef}\label{ndef1.2}
We denote by $A$ the largest possible constant such that
\[x^\frac{1}{37}\geq A\log x\quad\mbox{for all }x\geq 2.\]
\end{ndef}
Next we state a key lemma due to L. Pyber \cite{pyber} and an
immediate consequence which can also be found in
\cite{bertramsurvey}.\\

\begin{nlemma}\label{nlem1.3}
(a) Let $G$ be a group with $|G|\geq 4$. Suppose that $N\unlhd G$
is nilpotent and that
\[k(G/N)\geq 2^{x(\log|G/N|)^\frac{1}{t}}\]
for some constants $0<x\leq 1$ and $t\geq 1$. Then
\[k(G)\geq\frac{x^t \log|G|}{2(\log\log|G|)^t}.\]
(b) If $G$ is a group with $|G|\geq 4$ and $N\unlhd G$ is nilpotent such that
$k(G/N)\geq |G/N|^x$ for some $0<x\leq 1$, then
\[k(G)\geq\frac{x}{2}\frac{\log|G|}{\log\log|G|}.\]
\end{nlemma}
\begin{beweis}
For (a) see \cite[Lemma 2.2]{pyber}, and (b) is the special case of (a) where
$t=1$.
\end{beweis}

We also recall the following results:\\

\begin{nsatz}\label{nsatz1.4}
Let $0\not=V$ be a faithful, completely reducible, finite $G$-module for a solvable group $G$. Then
\[|G|\leq\frac{|V|^\alpha}{\lambda}\]
where $\lambda=24^\frac{1}{3}=2\cdot\sqrt[3]{3}$ and $9^\alpha=48\lambda$, so that
$2.24<\alpha<2.25$.
\end{nsatz}
\begin{beweis}
See \cite[Theorem 3.5(a)]{manz/wolf}.
\end{beweis}

\begin{nlemma}\label{nlem1.5}
There is an absolute constant $\alpha>0$ such that if $G$ is a solvable permutation group on a faithful set
$\Omega$, then
\[s(G)\geq|G|^\alpha.\]
\end{nlemma}
\begin{beweis}
This follows immediately from \cite[Lemma 2.4]{pyber}. (Note that our $\alpha$ here differs slightly from the
$\alpha$ in Pyber's paper.)
\end{beweis}

We remark that the previous lemma is a special case of a more general (but elementary) result of Babai
and Pyber \cite{babai-pyber} and that no explicit value of $\alpha$ is given in their
work.\\

\begin{nlemma}\label{nlem1.6}
Let $G\not=1$ be a group and $N$ a normal subgroup of $G$. Suppose that $\alpha,\beta$ are
real numbers with $0\leq\alpha\leq 1$ and $0\leq\beta\leq 1$ such that $k(N)\geq|N|^\alpha$ and
$k(G/N)\geq|G/N|^\beta$. Then
\[k(G)>|G|^\frac{\alpha\beta}{1+\alpha+\beta}\geq
|G|^\frac{\alpha\beta}{3}.\]
\end{nlemma}
\begin{beweis}
See \cite[Lemma 3(i)]{bertramisrael}.
\end{beweis}

\section{The main results}\label{section2}

We first prove an improvement of a result of Pyber (see \cite[Theorem B]{pyber}).
Pyber proved that if $G$ is solvable with $\Phi(G)=1$ and $|G|\geq 4$, then
$k(G)\geq 2^{\beta\log|G|/(\log\log|G|)^2}$ for some universal constant $\beta>0$. Here
we provide a polynomial lower bound in $|G|$ in the same
situation.\\

\begin{nsatz}\label{nsatz2.1}
Let $G$ be a solvable group with $\Phi(G)=1$. Then
\[k(G)\geq|G|^\beta\]
for some universal constant $\beta>0$. Specifically, one can choose
$\beta=\frac{\alpha}{1191(\alpha+1)+1}$ where $\alpha>0$ is as in \ref{nlem1.5}.
\end{nsatz}
\begin{beweis}
We will use the well-known structural facts on solvable groups listed in
\cite[Lemmas 3.1 and 3.7]{pyber} and also use some of the ideas from the proof of
\cite [Theorem B]{pyber}.\\

As is well-known, since $\Phi(G)=1$, $V=F(G)$ can be considered to
be a faithful, completely reducible $H$-module of mixed
characteristic, where $H$ is a complement of $F(G)$ in $G$. So we
write $G=HV$ and write $V$ additively. More precisely, we can write
$V=V_1\oplus\ldots\oplus V_n$ for some $n\in\n$ where each $V_i$ is
an irreducible $\GF (p_i)G$-module for a prime $p_i$
($i=1,\ldots,n$). (Note that the $p_i$ need not be mutually
distinct.) Put $G_i=G/C_G(V_i)\cong H/C_H(V_i)$ for all $i$ and
relabel the $V_i$ in such a way that $n(G_i,V_i)\geq
|V_i|^\frac{1}{37}$ for all $i=1,\ldots,r$ for some $r\in\{0, 1, \dots ,n\}$, and
$n(G_i,V_i)<|V_i|^\frac{1}{37}$ for $i=r+1, \ldots,n$. Note that by
\ref{nsatz1.1}, for $i=r+1,\ldots,n$, we know that
$G_i$ has a normal subgroup $N_i$ of derived length at most 2 such that
$G_i/N_i$ is isomorphic to a subgroup of $S_{k_i}$ for some $k_i\in\n$,
and we can write
\[(*)\quad V_i=V_{i1}\oplus\ldots\oplus V_{i, k_i}\]
with $N_i$-modules $V_{ij}$ which are permuted by $G_i/N_i$, so that
\[N_i=\bigcap_{j=1}^{k_i}N_{G_i}(V_{ij})\cong
\bigcap_{j=1}^{k_i}N_G(V_{ij})/C_G(V_i).\]
Now let $W_1=V_1\oplus\ldots\oplus V_r$ and $W_2=V_{r+1}\oplus\ldots\oplus V_n$.
Clearly $W_i\unlhd G$ ($i=1,2$). Put
\[H_1=C_H(W_2)\ \mbox{and }N=H_1W_1\unlhd G.\]
Obviously $H_1$ acts faithfully on $W_1$. Next we put $H_2=H/H_1$ so that
$H_2$ acts faithfully on $W_2$. Write $T=H_2W_2$ for the semidirect product
of $H_2$ and $W_2$ with respect to that action and observe that $G/N\cong
T$.(Note that if $r=0$, then $W_1=0$, $W_2=V$, $H_1=1$, $N=1$ and $T=G$; and if
$r=n$, then $W_1=V$, $W_2=0$, $H_1=H$, $N=G$ and $T=1$.)\\

In view of \ref{nlem1.6}, we now seek lower bounds for $k(T)$ and $k(N)$
separately. For $k(N)$, observe that clearly
\begin{eqnarray*}
k(N) &\geq& n(H_1,W_1)\ \geq\ \prod_{i=1}^r n(H_1,V_i)\ \geq\ \prod_{i=1}^r n(G,V_i)\\
     &=&\prod_{i=1}^r n(G_i,V_i)\ \geq\ \prod_{i=1}^r |V_i|^\frac{1}{37}\ =\ |W_1|^\frac{1}{37}.
\end{eqnarray*}
Moreover, by \ref{nsatz1.4} we have $|H_1|\leq|W_1|^3$, and thus
$|N|\leq |W_1|^4$, so altogether
\[k(N)\geq
|W_1|^\frac{1}{37}\geq\left(|N|^\frac{1}{4}\right)^\frac{1}{37}=|N|^\frac{1}{148}\quad (1).\]
Next we study $k(T)$ and obtain a lower bound similarly as in the proof of \cite[Theorem B]{pyber}.
If $i\in\{r+1,\ldots,n\}$, then we observe that
$G_i\cong H/C_H(V_i)$ and thus let $M_i\leq H$ be the inverse image of $N_i$ in $H$.
Hence $M_i/C_H(V_i)=N_i$ and hence
\[M_i=\bigcap_{j=1}^{k_i}N_H(V_{ij}).\]
Put
\[\Omega=\{V_{ij}\ |\ i=r+1,\ldots,n\ ;j=1,\ldots,k_i\}.\]
Clearly $H$ acts on $\Omega$, and if we let
\[K=\bigcap_{i=r+1}^n\bigcap_{j=1}^{k_i}N_H(V_{ij})=\bigcap_{i=r+1}^n
M_i\]
be the kernel of this action, then $H/K\ug\Sym(\Omega)$. We also claim that
\[(2)\quad K/H_1\ug \X_{i=r+1}^n N_i.\]
To see this, observe that clearly $1=\bigcap\limits_{i=1}^n C_{H/H_1}(V_i)$.
Therefore, if for $U\leq H$ we write $\overline{U}=UH_1/H_1$, then
\begin{eqnarray*}
\overline{K}&\cong&\overline{K}/\bigcap\limits_{i=r+1}^n(\overline{K}\cap C_{\overline{H}} (V_i))\\
            &\ug&\X_{i=r+1}^n \overline{K}/(\overline{K}\cap
                    C_{\overline{H}}(V_i))\\
            &\cong&\X_{i=r+1}^n (\overline{K}
            C_{\overline{H}}(V_i))/C_{\overline{H}}(V_i)\\
            &\leq&\X_{i=r+1}^n \overline{M_i}/C_{\overline{H}}(V_i)\\
            &=&\X_{i=r+1}^n (M_i/H_1)/(C_H(V_i)/H_1)\\
            &\cong&\X_{i=r+1}^n M_i/C_H(V_i)\\
            &=&\X_{i=r+1}^n N_i
\end{eqnarray*}
which proves the claim.\\
Now as $\dl(N_i)\leq 2$ for $i=r+1,\ldots,n$, from (2) we see that for the semidirect product $\overline{K}W_2\leq T$ we have
\[\dl(\overline{K}W_2)\leq\dl(\overline{K})+\dl(W_2)\leq 2+1=3.\]
Now a result of Bertram \cite[Theorem 1]{bertramisrael} yields
\[k(\overline{K}W_2)\geq|\overline{K}W_2|^\frac{1}{7}\quad (3)\]
On the other hand,
\[S:=T/(\overline{K}W_2)\cong\overline{H}/\overline{K}\cong
H/K\ug\Sym(\Omega)\]
and thus by \ref{nlem1.5} we have
\[s(S)\geq |S|^\alpha\]
with the $\alpha$ as in \ref{nlem1.5}. Now if two elements $x_k\in W_2$ ($k=1,2$) are conjugate in $G$,
then they are conjugate in $T$, and if we write
\[v_k=\sum_{i=r+1}^n\sum_{j=1}^{k_i}\alpha_{kij}v_{kij}\]
for suitable $v_{kij}\in V_{ij}$ and $\alpha_{kij}$ in the field belonging to $V_{ij}$ ($k=1,2$), and
$t\in T$ is such that $v_1^t=v_2$, then clearly for the subsets
$\Omega_k=\{V_{ij}|\alpha_{kij}v_{kij}\not=0\}$ of $\Omega$ ($k=1,2$) we have that $\Omega_1^{t\overline{K}W_2}=\Omega_2$.
Therefore
\[k(T)\geq n(\overline{H},W_2)\geq s(S)\geq|S|^\alpha.\quad (4)\]
Now with (3) and \cite[Lemma 2.1(ii)]{pyber} we also get
\[k(T)\geq\frac{|\overline{K}W_2|^\frac{1}{7}}{|S|}\quad (5),\]
and putting (4) and (5) together gives us
\[|T|=|S||\overline{K}W_2|\leq|S| |S|^7k(T)^7=|S|^8k(T)^7\leq
k(T)^\frac{8}{\alpha}k(T)^7=k(T)^{\frac{8}{\alpha}+7}.\]
Therefore
\[k(G/N)=k(T)\geq|T|^\frac{1}{\frac{8}{\alpha}+7}=|T|^\frac{\alpha}{8+7\alpha}\quad
(6).\]
Finally, using (1), (6), and \ref{nlem1.6}, we obtain
\[k(G)\geq|G|^\beta \ \mbox{with
}\beta=\frac{\alpha}{1148+1036\alpha}/\left(\frac{149}{148}+\frac{\alpha}{8+7\alpha}\right)=
\frac{\alpha}{1191(\alpha+1)+1}\ ,\]
as desired.
\end{beweis}

Our first application are improvements on the strongest general lower bounds for $k(G)$ in
terms of $|G|$ for solvable and arbitrary finite groups, as obtained by Pyber in
\cite{pyber}.\\
As to solvable groups, Pyber's arguments show that
\[k(G)\geq\beta_1\frac{\log|G|}{(\log\log|G|)^3}\]
for a suitable (universal) constant $\beta_1$ (see \cite[Corollary 2.1]{bertramsurvey}).
Here we improve this as follows.\\

\begin{ncor}\label{ncor2.2}
Let $G$ be a solvable group with $|G|\geq 4$, and let $\beta$ be as in \ref{nsatz2.1}. Then
\[k(G)\geq\frac{\beta}{2}\frac{\log|G|}{\log\log|G|}.\]
\end{ncor}
\begin{beweis}
Let $N=\Phi(G)$. By \ref{nsatz2.1} we know that $k(G/N)\geq |G/N|^\beta$,
so by \ref{nlem1.3}(b) we get the desired conclusion.
\end{beweis}

For arbitrary groups, the best possible bound to date is Pyber's result that there
exists an $\epsilon>0$ such that
\[k(G)\geq\epsilon\frac{\log|G|}{(\log\log|G|)^8}.\]
Using Theorem 2.1 in Pyber's proof, we can get the following slight
improvement:\\

\begin{ncor}\label{ncor2.3}
There exists a constant $\epsilon_1>0$ such that if $G$ is a group with $|G|\geq 4$, then
\[k(G)\geq\epsilon_1\frac{\log|G|}{(\log\log|G|)^7}.\]
\end{ncor}
\begin{beweis}
Let $O=O_\infty(G)$ be the largest normal solvable subgroup of $G$,
so clearly $O_\infty(G/O)=1$. Let $\Phi=\Phi(O)$, and let $x$ and $y$
denote the orders of the factor groups $O/\Phi$ and $G/O$, respectively. Then $\Phi(O/\Phi)=1$,
and thus by \ref{nsatz2.1} we have $k(O/\Phi)\geq x^\beta$,
where $\beta$ is as in \ref{nsatz2.1}. If follows that
\[k(G/\Phi)\geq\frac{k(O/\Phi)}{|G/O|}\geq\frac{x^\beta}{y}.\]
Thus
\[(1)\quad\log k(G/\Phi)\geq\beta\log x-\log
y=\beta\log(xy)-(1+\beta)\log y=\beta\log|G/\Phi|-(1+\beta)\log y.\]
Now by \cite[Lemma 4.7]{pyber} we have
\[(2)\quad \log k(G/\Phi)\geq\log k(G/O)\geq\delta(\log
y)^\frac{1}{7}.\]
Combining (1) and (2) therefore yields
\begin{eqnarray*}
\beta\log|G/\Phi|&\leq&\log k(G/\Phi)+(1+\beta)\log y\\
                 &\leq&\log k(G/\Phi)+(1+\beta)\left(\frac{
                       \log k(G/\Phi)}{\delta}\right)^7.
\end{eqnarray*}
This easily implies that
\[\log k(G/\Phi)\geq\gamma(\log|G/\Phi|)^\frac{1}{7}\]
for some $0<\gamma\leq 1$. Thus
\[k(G/\Phi)\geq 2^{\gamma(\log|G/\Phi|)^\frac{1}{7}},\]
and as $\Phi$ is nilpotent, we may apply \ref{nlem1.3}(a) which yields
the assertion.
\end{beweis}

\section{More applications for solvable groups}\label{section3}

In this section we significantly strengthen some results of Bertram
\cite{bertramsurvey} on solvable groups.\\
In \cite[Proposition 2.3]{bertramsurvey} it is shown that if $G$ is
solvable, $\Phi(G)$ is abelian and $|G|$ is sufficiently large
(depending only on $t$), then $k(G)>(\log |G|)^t$. We now strengthen
the lower bound and generalize the hypothesis to allow $\Phi(G)$ to be of
a fixed (arbitrary) nilpotency class. What we get somewhat resembles
Sherman's bound $k(P)\geq c|P|^\frac{1}{c}-c+1$ for nilpotent groups
of nilpotency class $c$ (see \cite{sherman}).

\begin{nsatz}\label{nsatz3.1}
Let $G$ be a solvable group, and suppose $\Phi(G)$ is of nilpotency
class $c$. Then
\[\log
k(G)>\frac{\beta\left(\frac{1}{c}+\frac{\log(c/2)}{\log|\Phi|}\right)}{
1+\beta+\frac{1}{c}+\frac{\log(c/2)}{\log|\Phi|}}\log|G|,\]
where $\beta$ is as in \ref{nsatz2.1}. In particular,
\[k(G)>\left(\frac{1}{2}c|G|^\frac{1}{2}\right)^\frac{\beta}{3}.\]
Moreover, if $\Phi(G)$ is abelian, then
\[k(G)>|G|^\frac{\beta}{2+\beta}\geq |G|^\frac{\beta}{3}.\]
\end{nsatz}
\begin{beweis}
Write $\Phi=\Phi(G)$. By \ref{nsatz2.1} we have $k(G/\Phi)\geq |G/\Phi|^\beta$.
Also, by Sherman's result \cite{sherman} we have  $k(\Phi)\geq c|\Phi|^\frac{1}{c}-c+1$.
As $|\Phi|^\frac{1}{c}\geq 2$, it is easy to see that this implies that
$k(\Phi)\geq\frac{c}{2}|\Phi|^\frac{1}{c}=|\Phi|^{\frac{1}{c}+\frac{\log\frac{c}{2}}{\log|\Phi|}}$.
Thus by \ref{nlem1.6} we conclude that
\[\log
k(G)>\frac{\beta\left(\frac{1}{c}+\frac{\log(c/2)}{\log|\Phi|}\right)}{
1+\beta+\frac{1}{c}+\frac{\log(c/2)}{\log|\Phi|}}\log|G|
\geq\frac{\frac{\beta}{c}+\frac{\beta\log(c/2)}{\log|\Phi|}}{3}\log|G|.\]
Hence
\[k(G)\geq|G|^\frac{\beta}{3c}\cdot 2^{\frac{\beta}{3}\log(c/2)}=
|G|^\frac{\beta}{3c}\left(\frac{c}{2}\right)^\frac{\beta}{3}=
\left(\frac{1}{2}c|G|^\frac{1}{c}\right)^\frac{\beta}{3},\]
as claimed. Finally, if $\Phi$ is abelian, then $k(\Phi)=|\Phi|$,
and then \ref{nlem1.6} yields
\[k(G)>|G|^\frac{\beta}{2+\beta}\geq
|G|^\frac{\beta}{3},
and the proof is complete.          \]
\end{beweis}

Our next goal is to strengthen \cite[Corollary 2.3]{bertramsurvey}.\\

\begin{ncor}\label{ncor3.2}
Let $G$ be solvable, and write $F=F(G)$ for the Fitting subgroup of $G$.
Let $\beta$ be as in \ref{nsatz2.1}.\\
(a) Let $0<\alpha\leq 1$. If $|F'|\leq|G|^{1-\frac{\alpha(2+\beta)}{\beta}}$,
then $k(G)>|G|^\alpha$.\\

(b) Let $t>0$. If $|F'|\leq \frac{|G|}{(\log|G|)^{t(1+\frac{2}{\beta})}}$,
then $k(G)>(\log|G|)^t$.
\end{ncor}
\begin{beweis}
As in the proof of \cite[Corollary 2.3]{bertramsurvey} we see that
\[\Phi(G/F')=\Phi(G)/F'<F/F',\]
and thus $\Phi(G/F')$ is abelian. By \ref{nsatz3.1} we conclude that $k(G)\geq k(G/F')>|G/F'|^\frac{
\beta}{2+\beta}$.
Now $|G/F'|^\frac{\beta}{2+\beta}\geq |G|^\alpha$ if and only if $|F'|\leq |G|^{1-\frac{\alpha(2+\beta)}{\beta}}$,
as can easily be checked. This is (a). \\
Likewise, $|G/F'|^\frac{\beta}{2+\beta}\geq(\log|G|)^t$ if and only if $|F'|$ satisfies the hypothesis in
(b), and so we are done.
\end{beweis}

We finally improve \cite[Corollary 2.4]{bertramsurvey}.\\

\begin{ncor}\label{ncor3.3}
Let $G$ be a solvable Frobenius group with Frobenius kernel $N$. If $N$ is abelian, then
\[k(G)>|G|^\frac{\beta}{2+\beta}.\]
where $\beta$ is as in \ref{nsatz2.1}.
\end{ncor}
\begin{beweis}
As in the proof of \cite[Corollary 2.4]{bertramsurvey} we note that $N=F(G)$, and so we are done by
\ref{nsatz3.1}.
\end{beweis}

As already remarked in \cite{bertramsurvey}, note that if $G$ is a Frobenius group with Frobenius
kernel $N$ such that $|G/N|$ is even, then $N$ necessarily is abelian and \ref{ncor3.3} can be applied.

\end{document}